\documentclass[12pt]{article}
\usepackage{amsmath}
\usepackage{amssymb}
\usepackage{amscd}
\usepackage{latexsym}
\usepackage{theorem}
\usepackage{doublespace}
\usepackage{epsfig}
\usepackage{psfrag}
\usepackage{amsfonts}
\usepackage{enumerate}

\newtheorem{theorem}{Theorem}[section]
\newtheorem{lemma}[theorem]{Lemma}
\newtheorem{proposition}[theorem]{Proposition}
\newtheorem{corollary}[theorem]{Corollary}
\newtheorem{claim}[theorem]{Claim}

{\theorembodyfont{\rmfamily}
\theoremstyle{plain}

\newtheorem{example}[theorem]{Example}
\newtheorem{problem}{Problem}
\newtheorem{remark}[theorem]{Remark}
}

\input xy
   \xyoption{all}

\newcommand{\sz}{|S|-1}

\renewcommand{\>}{\right>}

\newcommand{\kd}{\mathfrak{k}^*}

\newcommand{\tr}{\operatorname{tr}}

\newcommand{\Ad}{\operatorname{Ad}}

\renewcommand{\dim}{\operatorname{dim}}

\newcommand{\Pd}{\displaystyle\prod}

\newcommand{\C}{{\mathbb{C}}}
\newcommand{\Zt}{{\mathbb{Z}_2}}
\newcommand{\Z}{{\mathbb{Z}}}
\newcommand{\Q}{{\mathbb{Q}}}
\newcommand{\R}{{\mathbb{R}}}

\newcommand{\la}{\lambda}
\renewcommand{\a}{\alpha}

\renewcommand{\k}{\kappa}
\newcommand{\kso}{\k_{S^1}}

\newcommand{\ctn}{{\C^{2n}}}

\newcommand{\tnd}{(\mathfrak{u}(1)^*)^n}
\newcommand{\SO}{{SO(3)}}
\newcommand{\SL}{{SL(2,\C)}}

\newcommand{\sutd}{{\mathfrak{su}(2)}^*}
\newcommand{\sut}{{\mathfrak{su}(2)}}
\newcommand{\slt}{{\mathfrak{sl}(2,\C)}}
\newcommand{\Ae}{[A;e_1,\ldots,e_n]}

\newcommand{\polsa}{M_S}
\newcommand{\pola}{M}
\newcommand{\mcs}{\mathcal{S}}
\newcommand{\sa}{\mcs}
\newcommand{\spa}{\mcs'}
\newcommand{\core}{\mathfrak{L}}

\newcommand{\subs}{\subseteq}
\newcommand{\sups}{\supseteq}
\newcommand{\hookto}{{\hookrightarrow}}

\newcommand{\impl}{\Rightarrow}

\newcommand{\hs}{\hspace{3pt}}

\renewcommand{\ss}{\substack}

\newcommand{\half}{\frac{1}{2}}
\renewcommand{\mod}{{/\!\!/\!}}
\newcommand{\mmod}{{/\!\!/\!\!/\!\!/\!}}

\newcommand{\mc}{\mu_{\C}}

\newcommand{\hgso}{H^*_{K\times S^1}}
\newcommand{\hg}{H^*_K}
\newcommand{\hso}{H^*_{S^1}}
\renewcommand{\cot}{T^*\C^{2n}}
\renewcommand{\L}{\mathcal{L}}
\newcommand{\K}{\mathcal{K}}

\newcommand{\ns}{n_S}
\newcommand{\bi}{b_i}
\newcommand{\bk}{b_k}
\newcommand{\bj}{b_j}
\newcommand{\ba}{b_A}
\newcommand{\bns}{b_{\ns}}

\newcommand{\ck}{c_k}

\newcommand{\cns}{c_{\ns}}
\newcommand{\Sb}{\overline{S}}
\newcommand{\Lb}{\overline{S^c}}
\newcommand{\Ac}{A^c}
\newcommand{\asb}{A\cap\Sb}
\newcommand{\acsb}{\Ac\cap\Sb}
\newcommand{\pis}{\prod_{i\in\Sb}}
\newcommand{\pjl}{\prod_{j\in\Lb}}
\newcommand{\vs}{v_S}

\newcommand{\on}{\{1,\ldots,n\}}
\newcommand{\twn}{\{2,\ldots,n\}}

\newcommand{\us}{U_S}
\newcommand{\mtus}{X_T\cap\us}
\newcommand{\ts}{T\sups S}

\newcommand{\mur}{\mu_{\R}}
\newcommand{\muc}{\mu_{\C}}
\newcommand{\ses}{\a-\text{st}}
\newcommand{\bigmid}{\hs\hs\big{|}\hs\hs}
 
\newcommand{\qed}{\hfill \mbox{$\Box$}\medskip\newline}
\newenvironment{proof}{\noindent {\bf Proof:}}{\qed \par}
\newenvironment{componentproof}{\noindent {\bf Proof of \ref{component}:}}{\qed \par}
\newenvironment{usproof}{\noindent {\bf Proof of \ref{us}:}}{\qed \par}
\newenvironment{mainproof}{\noindent {\bf Proof of \ref{main}:}}{\qed \par}
\newenvironment{eqcoreproof}{\noindent {\bf Proof of \ref{eqcore}:}}{\qed \par}

\begin{document}
\begin{spacing}{1.1}

\noindent
{\LARGE \bf Hyperpolygon spaces and their cores}\bigskip\\
{\bf Megumi Harada} \\
Department of Mathematics, University of Toronto,
Ontario M5S 3G3 Canada.\smallskip \\
{\bf Nicholas Proudfoot } \\
Department of Mathematics, University of California,
Berkeley, CA 94720.
\bigskip
{\small
\begin{quote}
\noindent {\em Abstract.}
Given an $n$-tuple of positive real numbers $(\a_1,\ldots,\a_n)$,
Konno \cite{K2} defines the {\em hyperpolygon space} $X(\a)$, 
a hyperk\"ahler
analogue of the K\"ahler variety $M(\a)$ parametrizing polygons in $\R^3$
with edge lengths $(\a_1,\ldots,\a_n)$.  The polygon space $M(\a)$
can be interpreted as the moduli space of stable representations of a 
certain quiver with fixed dimension vector; 
from this point of view, $X(\a)$ is the hyperk\"ahler
quiver variety defined by Nakajima \cite{N1,N2}.
A quiver variety admits a natural $\C^*$-action, and the union
of the precompact orbits is called the 
{\em core}.  We study the components of the core of $X(\a)$,
interpreting each one as a moduli space of pairs of polygons in $\R^3$
with certain properties.
Konno gives a presentation of the cohomology ring of $X(\a)$;
we extend this result by computing the $\C^*$-equivariant cohomology
ring, as well as the ordinary and equivariant cohomology rings of the
core components.
\end{quote}
}
\bigskip

Let $K$ be a compact Lie group acting linearly on $\C^N$
with moment map $\mu:\C^N\to\kd$ such that $\mu(0)=0$.
Then for any central regular value $\a\in\kd$, the K\"ahler quotient
$$M(\a) = \C^N\mod_{\a}K = \mu^{-1}(\a)/K$$ is a K\"ahler manifold
of complex dimension $N-\dim K$.
The cotangent bundle $T^*\C^N$ is a hyperk\"ahler manifold,
and the induced action of $K$ on $T^*\C^N$ is hyperhamiltonian
(see, for example, \cite{HP}), with hyperk\"ahler moment map
$$\mur\oplus\muc:T^*\C^N\to\kd\oplus\kd_{\C}.$$
We call the hyperk\"ahler reduction
$$X(\a) = T^*\C^N\mmod_{(\a,0)}K = \Big(\mur^{-1}(\a)\cap\muc^{-1}(0)\Big)\Big/K$$
the {\em hyperk\"ahler analogue} of $X(\a)$.
The manifold $X(\a)$ is a noncompact 
hyperk\"ahler manifold of complex dimension $2(N-\dim K)$,
containing the cotangent bundle to $M(\a)$ as an open set \cite{HP}.
The action of the nonzero complex numbers $\C^*$ on $T^*\C^N$
given by scalar multiplication on each fiber induces an action of $\C^*$
on $X(\a)$, which restricts to the scalar action on the fibers of
$T^*M(\a)\subs X(\a)$.
The {\em core} $\mathfrak{L}$ of $X(\a)$ is defined to be the set of points
$x\in X(\a)$ such that the limit 
$\lim_{\la\to\infty}\la\cdot x$ exists (the opposite limit
$\lim_{\la\to 0}\la\cdot x$ always exists).
The action of $\C^*$ defines a deformation retraction of $X(\a)$ onto $\mathfrak{L}$.
If $M(\a)$ is nonempty and compact, then $\mathfrak{L}$
is simply the union of all those Bia{\l}yniski-Birula strata
whose closures are compact.

In \cite{HP} we studied the hyperk\"ahler
analogues of toric varieties, 
which arise when $K$ is abelian (see also \cite{BD,HS,K1}). 
In this case, $\mathfrak{L}$ is a union of toric varieties, one of which
is the original toric variety $M(\a)$.
Another important context in which K\"ahler reductions and their
hyperk\"ahler analogues arise is the study of varieties associated to quivers;
this includes the spaces that we will study in this paper.

Suppose given a quiver $Q$ (a directed graph)
with vertex set $I$, and let $\{V_i\mid i\in I\}$ be a collection of 
finite complex dimensional vector
spaces.  A {\em representation} of $Q$ is a collection of maps
from $V_i$ to $V_j$ for every pair of vertices $i$ and $j$ connected by an edge.
The group $PU(V) = \Big(\prod GL(V_i)\Big)\Big/GL(1)_{\Delta}$ acts 
hamiltonianly on the space $E(Q,V)$
of representations of $Q$.
Both the K\"ahler quotients $M(\a) = E(Q,V)\mod_{\a}PU(V)$
and their hyperk\"ahler analogues $X(\a) = T^*E(Q,V)\mmod_{(\a,0)}PU(V)$, 
called {\em quiver varieties},
have been studied extensively.
A good introduction to these varieties, both
the K\"ahler and hyperk\"ahler versions, can be found in \cite{N2}.
Cores of quiver varieties have attracted particular attention in representation
theory.  The fundamental classes of their components provide a natural basis
for the top homology of $X(\a)$, which leads to the construction of canonical bases
for representations of the modified universal enveloping algebra associated
to the quiver $Q$ \cite{N3} (see also \cite{N4} for more recent results
along these lines).

\begin{figure}[h]
\begin{center}
\psfrag{a1}{$\a_1$}
\psfrag{a2}{$\a_2$}
\psfrag{a3}{$\a_3$}
\psfrag{a4}{$\a_4$}
\psfrag{a5}{$\a_5$}
\psfrag{a6}{$\a_6$}
\psfrag{a7}{$\a_7$}
\psfrag{a8}{$\a_8$}
\psfrag{a9}{$\a_9$}
\psfrag{0}{$\vec{0}$}
\includegraphics[height=40mm]{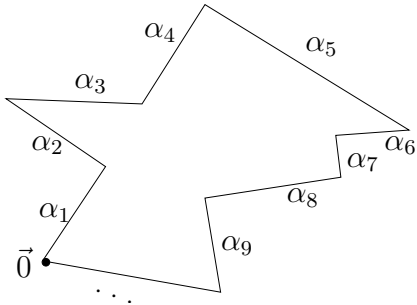}
\caption{A polygon in $\R^3$ with specified edge lengths $\a_i$.}\label{fig:polyex}
\end{center}\end{figure}

The examples with which we will be concerned in this paper are quiver varieties
corresponding to a very special class of quivers, as shown in the following picture.
\begin{figure}[h]
\centerline{\epsfig{figure=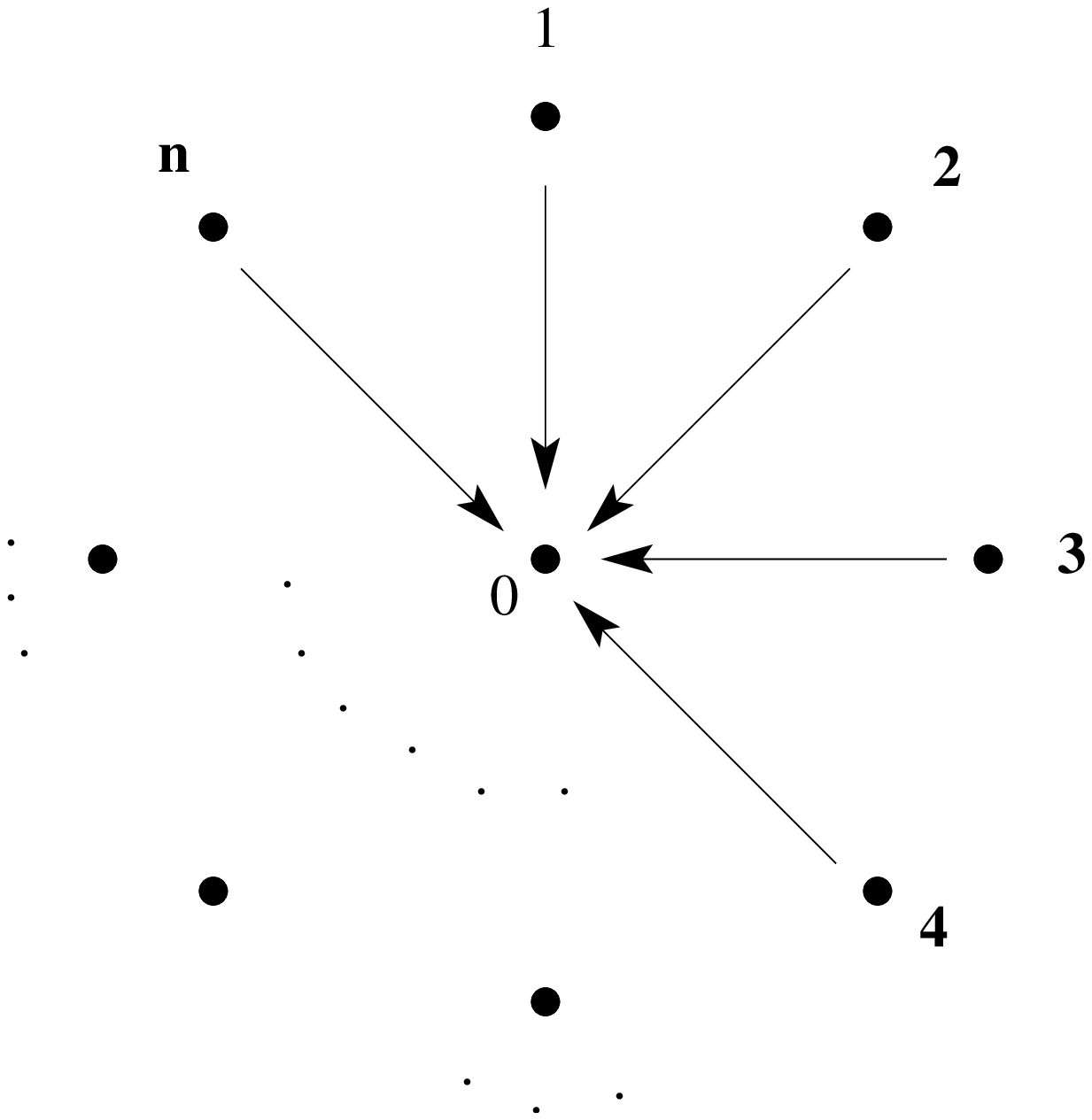,height=4cm}}
\caption{The quiver for hyperpolygon spaces.}\label{fig:quiver}
\end{figure}
We label the vertices $0$ through $n$, with zero in the center,
and for each $i\in\on$ we have an arrow from $i$ to $0$.
We put $V_i=\C$ for all $i\in\{1,\ldots,n\}$, and 
$V_0 = \C^2$, so that $$E(Q,V) = \bigoplus_{i=1}^n
\operatorname{Hom}(V_i,V_0)\cong\C^{2n}.$$
In this example, the K\"ahler quiver variety $M(\a)$
has a nice geometric interpretation.  The group
$PU(V)$ is isomorphic to $\Big(SU(2)\times U(1)^n\Big)\Big/\Z_2$
where $\Z_2$ acts diagonally on the $n+1$ factors,
and a central element $\a= 0 \oplus (\a_1,\ldots,\a_n)\in\sutd\oplus\tnd$
is given by an $n$-tuple of real numbers.
The variety $M(\a)$ is diffeomorphic to the moduli space of $n$-sided polygons in 
$\R^3$, with edge lengths $(\a_1,\ldots,\a_n)$, modulo the action of $SO(3)$
on $\R^3$ by rotation \cite{HK1, HK2, Kl}, as in Figure \ref{fig:polyex}.

The hyperk\"ahler analogue $X(\a)$ was introduced in \cite{K1},
in which Konno enumerated the components of the core
of $X(\a)$ (Theorem \ref{fixed}), showing that the components other than $M(\a)$
are in bijection with the collection of subsets $S\subs\on$ of size at least $2$
such that $\sum_{i\in S}\a_i < \sum_{j\in S^c}\a_j$.  Such a subset $S$ will be called
{\em short}.
Our first set of results, comprising Section \ref{core},
concerns the core components $\{U_S\}$, which one may think of as
{\em generalized polygon spaces}.
We prove that $U_S$ is smooth for each short subset $S$ (Theorem \ref{component}),
and interpret it as the moduli space of pairs of polygons in $\R^3$
with certain geometric properties (Theorem \ref{us}).
This is therefore a solution, in the special case
of polygon spaces, to the following general problem.

\begin{problem}
Given any moduli
space $M$ that can be constructed as a K\"ahler reduction
(or GIT quotient) of complex affine space, is it possible
to interpret the core components of the hyperk\"ahler analogue
$X$ as moduli spaces corresponding to other, related moduli problems?
\end{problem}

Our second set of results, comprising Sections \ref{s1} and
\ref{last}, concerns the $S^1$-equivariant cohomology of $X(\a)$ and
$U_S$, where $S^1\subs\C^*$ is the unit circle.\footnote{The
$S^1$-equivariant cohomology is identical to the $\C^*$-equivariant
cohomology; we use $S^1$ in Section~\ref{s1} because at times it is
convenient to work with the real ADHM description of $X(\a)$, in which
context the full $\C^*$-action is difficult to write down explicitly.}
Konno computes the cohomology of $X(\a)$, showing in particular that
the hyperk\"ahler Kirwan map
$$\k:H_K^*(T^*\C^{2n})\to H_K^*\Big(\mur^{-1}(\a)\cap\muc^{-1}(0)\Big)
\cong H^*(X(\a))$$
in the case of hyperpolygon spaces is surjective (Theorem \ref{Konno}).
We generalize this result by computing the kernel of the {\em equivariant}
Kirwan map
$$\kso:H_{K\times S^1}^*(T^*\C^{2n})\to H_{K\times S^1}^*
\Big(\mur^{-1}(\a)\cap\muc^{-1}(0)\Big)
\cong H_{S^1}^*(X(\a))$$
(Theorem \ref{main}), which is also
surjective by Corollary \ref{surjective}.
In Section \ref{last}, we compute
the ordinary and equivariant cohomology of a core component $U_S$
(Theorem \ref{eqcore} and Corollary \ref{ordcore}).
All cohomology rings are taken with coefficients in $\Q$.

In general, the compactness of the fixed point set $X(\a)^{S^1}$
gives equivariant cohomology many nice properties not shared
by the ordinary cohomology of $X(\a)$; for example, the localization theorem 
of \cite{AB} makes possible a theory of integration in equivariant
cohomology, provided that the fixed point set is compact.
Nakajima studies the $S^1$-equivariant cohomology and $K$-theory
of quiver varieties in \cite{N4}, using it to construct representations
of quantum affine algebras associated to $Q$.
In the process, he conjectures that the equivariant Kirwan map
is surjective for all quiver varieties \cite[7.5.1]{N4}.

The reader is advised that cohomology computations
in hyperpolygon spaces and generalized polygon spaces are often
motivated by beautiful and intuitive geometry, but they are
just as often driven by daunting, labyrinthine algebra.
Whenever possible, we precede the proofs of our theorems with
remarks that are aimed at making the geometry maximally transparent.
\newline

\noindent {\em Acknowledgments.}
We are grateful to Allen Knutson for guiding 
us through the polygonal world,
and to Tam\'as Hausel for introducing us to the work of Konno.
Also to Michael Thaddeus
for useful discussions.
The second author thanks the Swiss National Funds for 
Scientific Research for its support.

\begin{section}{Hyperpolygon spaces}\label{hp}
We begin by collecting the basic definition and properties
of a hyperpolygon space, most of which can be found in \cite{K2}.
Fix a positive integer $n\geq 3$, and consider the group\footnote{One 
may prefer to just consider
the group $\SL\times (\C^*)^n$, and allow it to act with a finite kernel.
We quotient by $\Zt$ only to be consistent with the conventions of \cite{N2}
and \cite{K2}.}
$$G := \Big(\SL\times (\C^*)^n\Big)\Big/\Z_2,$$
where $\Zt$ acts by multiplying each factor by $-1$.
We define a right action of $G$ on $\ctn$ as follows.
We will write an element of $\ctn$ as an $n$-tuple
of column vectors
$$q = (q_1,\ldots,q_n),$$
and put
$$q[A;e_1,\ldots,e_n] = (A^{-1}q_1e_1,\ldots,A^{-1}q_ne_n).$$
The compact subgroup
$$K := \Big(SU(2)\times U(1)^n\Big)\Big/\Z_2 \subs G$$
acts with moment map $\mu:\ctn\to\sutd\oplus\tnd$
given by the equation $$\mu(q) = \sum_{i=1}^n (q_iq_i^*)_0 \oplus 
\left(\half|q_1|^2,\ldots,\half|q_n|^2\right),$$
where $q_i^*$ denotes the conjugate transpose of $q_i$, $(q_iq_i^*)_0$ 
denotes the
traceless part of $q_iq_i^*$, and $\sutd$ is identified with 
$i\cdot\sut$ via the trace form.
Given an $n$-tuple of real numbers $(\a_1,\ldots,\a_n)$, we define the {\em polygon space}
$$M(\a) := \ctn\mod_{\a}K = \mu^{-1}(\a)/K,$$
where $\a = 0\oplus (\a_1,\ldots,\a_n)\in\sutd\oplus\tnd$.
If we break the reduction into two steps, reducing first by $U(1)^n$
and then by $SU(2)$, we find that
\begin{equation}\label{polygon}
M(\a)\cong\left\{(v_1,\ldots,v_n)\in (\R^3)^n\,\bigg|\, 
\|v_i\|=\a_i\text{ and }\sum v_i=0\right\}\bigg/SO(3)
\end{equation}
(see Remark \ref{interp} and the proof of Theorem \ref{us}).
Here $\sutd$ is being identified with $\R^3$, and the coadjoint action of $SU(2)$
on $\sutd$ is being replaced by the standard action of $SO(3)$ on $\R^3$ \cite{HK2}.
This space, therefore, may be thought of as the moduli space of $n$-sided polygons in $\R^3$,
with fixed edge lengths, up to rotation.  
In particular, $M(\a)$ is empty unless $\a_i\geq 0$ for all $i$.

We call $\a$ {\em generic} if there does not exist a 
subset $S\subs\{1,\ldots,n\}$
such that $\sum_{i\in S}\a_i=\sum_{j\in S^c}\a_j$.
Geometrically, this means that there is no element of $M(\a)$ represented by a polygon that is
contained in a single line in $\R^3$.  If $\a$ is generic, then $M(\a)$ is smooth \cite{HK1}.
Throughout this paper we will assume that $\a$ is generic, and that $\a_i> 0$ for all $i$.

To define the hyperk\"ahler analogue of $M(\a)$, we consider the induced action of $G$
on $\cot$.
Explicitly, we write an element of $\cot$ as $(p,q)$,
where $q = (q_1,\ldots,q_n)$ is an $n$-tuple
of column vectors
and $p = (p_1,\ldots,p_n)$ an $n$-tuple
of row vectors, and we put
$$(p,q)[A;e_1,\ldots,e_n] = 
\big((e_1^{-1}p_1A,\ldots,e_n^{-1}p_nA),(A^{-1}q_1e_1,\ldots,A^{-1}q_ne_n)\big).$$
The vector space $\cot$ has the structure of a hyperk\"ahler manifold,
and the action of $K$ on $\cot$ is hyperhamiltonian
with hyperk\"ahler moment map \cite{K2} (see also \cite{HP})
$$\mur\oplus\muc:\cot\to\Big(\sutd\oplus\tnd\Big)\oplus\Big(\slt^*\oplus 
(\mathfrak{u}(1)^n_{\C})^*\Big)$$
given by the equations
\begin{equation*}\label{mur}
\mu_{\R}(p,q)  =  \frac{\sqrt{-1}}{2}\sum_{i=1}^n \left(q_i q_i^* - p_i^* p_i\right)_0
\oplus\left(\half\left(|q_1|^2 - |p_1|^2\right), \ldots, \half\left(|q_n|^2 - |p_n|^2\right)\right) 
\end{equation*}
and
\begin{equation*}\label{muc}
\mu_{\C}\left(p,q\right)  =  - \sum_{i=1}^n \left(q_i p_i\right)_0\oplus \left(\sqrt{-1} p_1 q_1,
\ldots, \sqrt{-1} p_n q_n\right).
\end{equation*}
We then define the {\em hyperpolygon space} to be the hyperk\"ahler quotient
$$X(\a) := \cot\mmod_{(\a,0)}K = \Big(\mur^{-1}(\a)\cap\muc^{-1}(0)\Big)\Big/K,$$
a smooth, noncompact hyperk\"ahler manifold of complex dimension $2(n-3)$ \cite{K2}.

The polygon and hyperpolygon spaces $M(\a)$ and $X(\a)$ are precisely the K\"ahler
and hyperk\"ahler varieties associated by Nakajima 
to the quiver shown in Figure~\ref{fig:quiver}.
It is shown in \cite{N1} that
$$M(\a)\cong\left(\ctn\right)^{\ses}/G\hspace{.7cm}\text{ and }\hspace{.7cm}
X(\a)\cong\muc^{-1}(0)^{\ses}/G,$$
where $\a$-st means stable with respect to the weight $\a$ in the sense of
geometric invariant theory.
Nakajima gives a stability criterion for general quiver varieties \cite{N1,N2},
which Konno interprets in the special case of hyperpolygon spaces.
Call a subset $S\subs\{1,\ldots,n\}$ {\em short} if $\sum_{i\in S}\a_i < \sum_{j\in S^c}\a_j$,
otherwise call it {\em long}.  (Assuming that $\a$ is generic is equivalent to assuming
that every subset is either short or long.)
Given a point $(p,q)\in\cot$ and a subset $S\subs\{1,\ldots,n\}$,
we will say that $S$ is {\em straight} 
in $(p,q)$ if $q_i$ is proportional to $q_j$ for every $i,j\in S$.
The terminology comes from K\"ahler polygon spaces, in which this condition
is equivalent to asking that the vectors $v_i$ and $v_j$ be proportional
over $\R_+$, or that the edges of lengths $\a_i$ and $\a_j$ (if they happen to be adjacent)
line up to make a single edge of length $\a_i+\a_j$, 
as in Figure \ref{straight}.
\begin{figure}[h]
\begin{center}
\psfrag{a1}{$\a_1$}
\psfrag{a2}{$\a_2$}
\psfrag{a3}{$\a_3$}
\psfrag{a4}{$\a_4$}
\psfrag{a5}{$\a_5$}
\psfrag{a6}{$\a_6$}
\psfrag{a7}{$\a_7$}
\psfrag{a8}{$\a_8$}
\includegraphics[height=40mm]{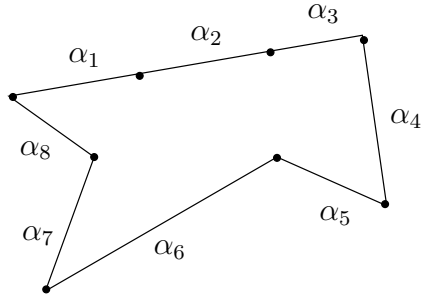}
\caption{The subset $\{1,2,3\}$ is straight.}\label{straight}
\end{center}\end{figure}

\begin{theorem}\label{stability}{\em\cite[4.2]{K2}}
Suppose that $\a$ is generic, and $\a_i> 0$ for all $i$.
Then a point $(p,q)\in\cot$ is $\a$-stable
if and only if the following two conditions are satisfied:
\begin{eqnarray*}
&1)& q_i\neq 0 \text{ for all }i,\text{ and}\\
&2)& \text{if $S$ is straight and $p_j=0$ for all $j\in S^c$, then $S$ is short.}
\end{eqnarray*}
\end{theorem}
We will use the notation $[p,q]$ to denote
the $G$-equivalence class of a point $(p,q)\in\muc^{-1}(0)^{\ses}$,
and $[p,q]_{\R}$ to denote
the $K$-equivalence class of a point $(p,q)\in\mur^{-1}(\a)\cap\muc^{-1}(0)$.
Note that $M(\a)$ sits inside of $X(\a)$ as the locus of points $[p,q]$ with $p=0$.
This observation, along with Theorem \ref{stability}, allows us to recover
the $\a$-stability condition for the action of $G$ on $\ctn$.
A point $q\in\ctn$ is $\a$-stable if and only if $q_i\neq 0$ for all $i$,
and no long subset is straight, as first shown in \cite{Kl}.
The polygonally-minded reader is warned that in the hyperpolygon space $X(\a)$,
long subsets {\em can} be straight.
\end{section}

\begin{section}{The core}\label{core}
For the rest of the paper we fix a generic $\a = 0\oplus (\a_1,\ldots,\a_n)\in\sutd\oplus\tnd$,
with $\a_i>0$ for all $i$.  To simplify notation, we will write $X=X(\a)$ and $M=M(\a)$.
Consider the action of $\C^*$ on $X$ given in the complex description
by $\la\cdot[p,q]=[\la p, q]$.  
The circle $S^1\subs\C^*$ preserves
$\omega_{\R}$, and acts with moment map $\Phi:X\to\R$ given in the symplectic
description by $\Phi\left([p,q]_{\R}\right) = \half\sum |p_i|^2$.
Following Konno, we define
$$\mathcal{S} = \big\{S\subs\{1,\ldots,n\}\hs\hs\big{|}\hs\hs S\text{ is short}\big\}$$
and $$\spa = \big\{S\in\mathcal{S}\hs\hs\big{|}\hs\hs |S|\geq 2\big\}.$$

\begin{theorem}\label{fixed}{\em\cite{K2}}
The fixed point set 
$\displaystyle{X^{\C^*} = X^{S^1} = \pola\cup\bigcup_{S\in\spa}X_S,}$
where
$$X_S = \big\{[p,q]\hs\hs\big{|}\hs\hs S\text{ and }S^c\text{ are each straight, and }
p_j = 0\text{ for all }j\in S^c\big\}.$$  Furthermore, $X_S$
is diffeomorphic to $\C P^{|S|-2}$.
\end{theorem}

For all $S\in\spa$, let
$U_S$ be the closure inside of $X$ of the set
$$\big\{[p,q]\hs\hs\big{|}\hs\hs\lim_{\la\to\infty}\la\cdot[p,q]\in X_S\big\},$$
and let $$\core = \pola\cup\bigcup_{S\in\spa}U_S.$$  This reducible subvariety
is called the {\em core} of $X$.  Since $\displaystyle{\lim_{\la\to 0}\la\cdot[p,q]}$ always exists,
the core is simply the union of those
Bia{\l}ynicki-Birula strata whose closures are compact.\footnote{In Morse theoretic language,
$U_S$ is the closed flow-down set for $X_S$ with respect to the Morse-Bott function $\Phi$,
and $C$ is the union of all of the flow-down sets.}
The $\C^*$ action defines an $S^1$-equivariant deformation retraction of $X$ onto $\core$ \cite{N1}.

\begin{theorem}\label{component}
The core component $U_S$ is smooth of complex
dimension $n-3$, and we have
$$U_S = \big\{[p,q]\bigmid S\text{ is straight, and }
p_j = 0\text{ for all }j\in S^c\big\}.$$
\end{theorem}

Before proving Theorem \ref{component}, we describe the way in which the various components
of the core fit together.
For all $S\in\spa$, let $$\polsa = U_S\cap \pola = 
\big\{[0,q]\hs\hs\big{|}\hs\hs S\text{ is straight}\big\}.$$
We call this space the {\em polygon subspace} of $\pola$ corresponding to the short subset $S$.
Note that $\polsa$ is itself a polygon space with $n-|S|+1$ edges, of lengths
$\{\a_j\mid j\in S^c\}\cup\{\sum_S\a_i\}$.  In particular, it is smooth.
Now suppose given two short subsets $S,T\in\spa$, and consider the intersection $U_S\cap U_T$.
\begin{itemize}
\item If $S\cap T=\emptyset$, then $U_S\cap U_T = \polsa\cap M_T$, a polygon subspace
both of $\polsa$ and of $M_T$.  
\item If $S\cap T\neq\emptyset$
and $S\cup T$ is long, then $U_S\cap U_T=\emptyset$.  
\item If $S\cap T\neq\emptyset$ and 
$S\cup T$ is short,
then $$U_S\cap U_T = \big\{[p,q]\bigmid S\cup T\text{ is straight, and }p_j=0\text{ for all }
j\in (S\cap T)^c\big\}.$$
This is a subvariety of $U_{S\cup T}$ given by taking the closure inside of $U_{S\cup T}$
of a certain subbundle of the conormal bundle to $M_{S\cup T}\subs\pola$,
defined by setting $p_j=0$ for all $j\in (S\cap T)^c\supseteq (S\cup T)^c$.
\end{itemize}
Each of these descriptions generalizes to higher intersections
without modification.

Finally, we compute the fixed point set $U_S^{\C^*}$.
If $[p,q]\in U_S^{\C^*}$, then either $p=0$ and $[p,q]\in\polsa$,
or $[p,q]\in X_T$ for some $T\in\spa$.  If $[p,q]\in X_T$ then Theorem \ref{fixed}
tells us that $T$ and $T^c$ are each straight, hence $S\subs T$ or $S\subs T^c$.
Since $p\neq 0$, we must have $S\subs T$.  Indeed, $U_S\cap X_T$ is the linear subspace
of $X_T\cong \C P^{|T|-2}$ given by the condition $p_j=0$ for all $j\in T\smallsetminus S$.
In particular, $U_S\cap X_T$ is isomorphic to $\C P^{|S|-2}$ for any $T\supseteq S$.

\begin{example}\label{ooh}
Let $n=5$, $\a_1=\a_2=1$, and $\a_3=\a_4=\a_5=3$, and consider the short subset
$S=\{1,2\}$.  The fixed point set of $U_S$ consists of $\polsa\cong\C P^1$,
and four points $X_S$, $U_S\cap X_{T_3}$, $U_S\cap X_{T_4}$, and $U_S\cap X_{T_5}$,
where $T_j = \{1,2,j\}$ for $j=3,4,5.$  
For each $j$, $U_S\cap U_{T_j}$ is isomorphic to $\C P^1$,
and touches $\polsa$ at the point $M_{T_j}$.
In the following picture, an ellipse represents a copy of $\C P^1$ flowing between two
fixed points, where the numbers or pairs of numbers
indicate subsets that are straight on this $\C P^1$.  (For example, 
$12,45$ means that $1$ and $2$ are straight,
as are $4$ and $5$.)
We will revisit this example at the end of Section \ref{last}.
\begin{figure}[h]
\begin{center}
\psfrag{Xs}{$X_S$}
\psfrag{Xt3}{$X_{T_3}$}
\psfrag{Xt4}{$X_{T_4}$}
\psfrag{Xt5}{$X_{T_5}$}
\psfrag{Ms(a)}{$M_{S}$}
\psfrag{Mt3(a)}{$M_{T_3}$}
\psfrag{Mt4(a)}{$M_{T_4}$}
\psfrag{Mt5(a)}{$M_{T_5}$}
\includegraphics[height=40mm]{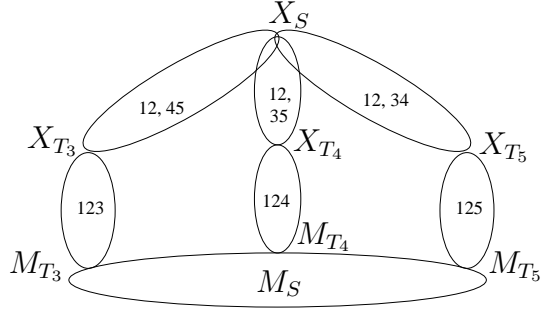}
\caption{$U_S$, with $S=\{1,2\}$}
\label{surface}
\end{center}
\end{figure}
\end{example}

\begin{componentproof}
The fact that $\dim U_S = \half\dim X = n-3$ is a general property of
core components of quiver varieties \cite{N1}.  
Thus, by a dimension count, it is enough to show that the set
$\big\{[p,q]\bigmid S\text{ is straight, and }
p_j = 0\text{ for all }j\in S^c\big\}$ is contained in $U_S$.

Consider a point $[p,q]\in X$ with $S$ straight, and $p_j=0$ for all $j\in S^c$.
By applying an element of $G$, we may assume that $q_i = \binom{1}{0}$ for all $i\in S$.
Suppose further that there exists an $i\in S$ with $p_i\neq 0$, and that
no strict superset of $S$ is straight.  In other words, if 
$q_j = \binom{a_j}{b_j}$ for $j\in S^c$,
suppose that $b_j\neq 0$.
For $\la\in\C^*$, let $A(\la) = 
\footnotesize{\Big(\begin{array}{cc}
\la&0\\
0&\la^{-1}
\end{array}
\Big)}$, 
let $e_i(\la)=\la$ for all $i\in S$,
and let $e_j(\la)=\la^{-1}$ for all $j\in S^c$.
Then for $i\in S$, we have $e_i(\la)^{-1}p_i A(\la) = \la^{-2} p_i$ and 
$A(\la)^{-1}q_i e_i = q_i$.
For $j\in S^c$, we have $A(\la)^{-1}q_j e_j = \binom{\la^{-2}a_j}{b_j}$.
Hence
\begin{eqnarray*}
\lim_{\la\to\infty}\la\cdot[p,q] &=& \lim_{\la\to\infty}\la^2\cdot[p,q]\\
&=& \lim_{\la\to\infty}[\la^2 p,q]\\
&=& \lim_{\la\to\infty}[\la^2 e(\la)^{-1}p A(\la), A(\la)^{-1}q e(\la)^{-1}]\\
&=& [p,q'],
\end{eqnarray*}
where $q'_i = q_i$ for $i\in S$, and $q'_j = \binom{0}{b_j}$ for $j\in S^c$.
Since we have assumed that $b_j\neq 0$ for all $j\in S^c$ and that $p_i\neq 0$
for some $i\in S$, $(p,q')$ is stable, and hence defines an element of $X_S$.
Since $U_S$ is defined to be the closure of the set of elements that flow up to $X_S$,
it includes all $[p,q]$ with $S$ straight and $p_j=0$ for all $j\in S^c$.

To see that $U_S$ is smooth, it is sufficient to show that $U_S$ is smooth
at $[p,q]$ for all $[p,q]\in X^{\C^*}$.
First suppose that $[p,q]\in X_T$ for some $T\in\spa$ containing $S$.
Suppose, without loss of generality, that
$T=\{1,\ldots,l\}$ and $S=\{1,\ldots,m\}$ for some $l\leq m$.
Konno computes an explicit local complex chart for $X$ at the point $[p,q]$,
with coordinates $\{z_i, w_i\mid 3\leq i\leq n-1\}$ \cite{K2}.
With respect to these coordinates, a point $[p',q']$ has the property that
$S$ is straight and $p'_j=0$ for all $j\in S^c$ if and only if 
$w_i=0$ for all $3\leq i\leq l$
and $z_j=0$ for all $l+1 \leq j\leq n-1$.
Hence $U_S$ is smooth at $[p,q]$.

It remains only to show that $U_S$ is smooth at $\polsa = U_S\cap\pola$.
Let $$E = \{(p,q)\mid S\text{ is straight, }
p_j = 0\text{ for all }j\in S^c,\text{ and }\mu_{\C}(p,q)=0\},$$ and let 
$N = \{(p,q)\in E\mid p=0\}$.  The natural projection from $E$ to $N$
exhibits $E$ as a vector bundle over $N$, because the equation $\mu_{\C}(p,q)=0$
is linear in $p$.  By definition, $U_S = E\mod G = E^{\ses}/G$, and 
$\polsa = N\mod G = N^{\ses}/G$.  
The set $E\vert_{N^{\ses}}/G\subs E^{\ses}/G$ is
an open neighborhood of $\polsa$ inside of $U_S$, and is isomorphic to a vector
bundle over $\polsa$.
Since $\polsa$ is a polygon space it is smooth, hence $U_S$ is smooth 
in a neighborhood of $\polsa$.
\end{componentproof}

The following corollary is known to the experts; 
we include it here for lack of an explicit reference.

\begin{corollary}\label{compactification}
$U_S$ is a compactification of the conormal
bundle to $\polsa$ in $\pola$.
\end{corollary}

\begin{proof}
We must show that the normal bundle to $\polsa$ in $U_S$
is dual to the normal bundle to $\polsa$ in $\pola$.
We use only general facts about quiver varieties from \cite{N1},
and the additional information that $U_S$ is smooth, 
from Theorem \ref{component}.
Consider a point $[0,q]\in\polsa$, and let $H_0$ and $H_1$ be the
$0$ and $1$ weight spaces of the $\C^*$ action on $T_{[0,q]}X$.  
The holomorphic symplectic
form $\omega_{\C}$ is being rotated by $\C^*$ with weight $1$ \cite[5.1]{N1},
hence it defines a perfect pairing between $H_0$ and $H_1$.
The fiber at $[0,q]$ of the normal bundle to $\polsa$ in $\pola$
is $H_0/T_{[0,q]}\polsa$, which is dual by $\omega_{\C}$ to the annihilator
of $T_{[0,q]}\polsa$ inside of $H_1$.
Since $U_S$ is $\C^*$-invariant, we may write 
$$T_{[0,q]}U_S = T_{[0,q]}U_S\cap H_0\oplus T_{[0,q]}U_S\cap H_1
= T_{[0,q]}\polsa\oplus T_{[0,q]}U_S\cap H_1.$$
To prove Corollary \ref{compactification}, we must show that
$T_{[0,q]}U_S\cap H_1$ is equal to the annihilator of $T_{[0,q]}\polsa$.
The fact that $T_{[0,q]}U_S\cap H_1$
is contained in the annihilator of $T_{[0,q]}\polsa$ follows from the
fact that $U_S$ is lagrangian with respect to $\omega_{\C}$ \cite{N1} 
(here we use
smoothness of $U_S$ at $[0,q]$).
Equality is then obtained by dimension count.
\end{proof}

We next describe $U_S$ in polygon-theoretic terms, as a certain moduli
space of pairs of polygons in $\R^3$.

\begin{theorem}\label{us}
Let $U_S$ be the component of the core of $X$ corresponding to a
subset $S\in\spa$.  Then $U_S$ is 
homeomorphic to the moduli space
of $n+1$ vectors $$\{u_i, v_j, w\in\R^3\mid i\in S, j\in S^c\},$$
taken up to rotation, satisfying the following conditions:
\begin{eqnarray*}
&1)&\hs\hs w + \sum_{j\in S^c}v_j = 0\\
&2)&\hs\hs \sum_{i\in S}u_i = 0\\
&3)&\hs\hs u_i\cdot w = 0\hs\hs\text{ for all }i\in S\\
&4)&\hs\hs \|v_j\| = \a_j\hs\hs\text{ for all }j\in S^c\\
&5)&\hs\hs \|w\| = \sum_{i\in S}\sqrt{\a_i^2 + \|u_i\|^2}.
\end{eqnarray*}
\end{theorem}

\begin{remark}\label{interp}
In more descriptive terms, a point in $U_S$ specifies two polygons in $\R^3$,
as in Figure~\ref{fig:Uspoly}.
\begin{figure}[h]
\begin{center}
\psfrag{w}{$w$}
\psfrag{v1}{$v_1$}
\psfrag{v2}{$v_2$}
\psfrag{v3}{$v_3$}
\psfrag{vnS}{$v_{n-|S|}$}
\psfrag{u1}{$u_1$}
\psfrag{u2}{$u_2$}
\psfrag{uS}{$u_{|S|}$}
\includegraphics{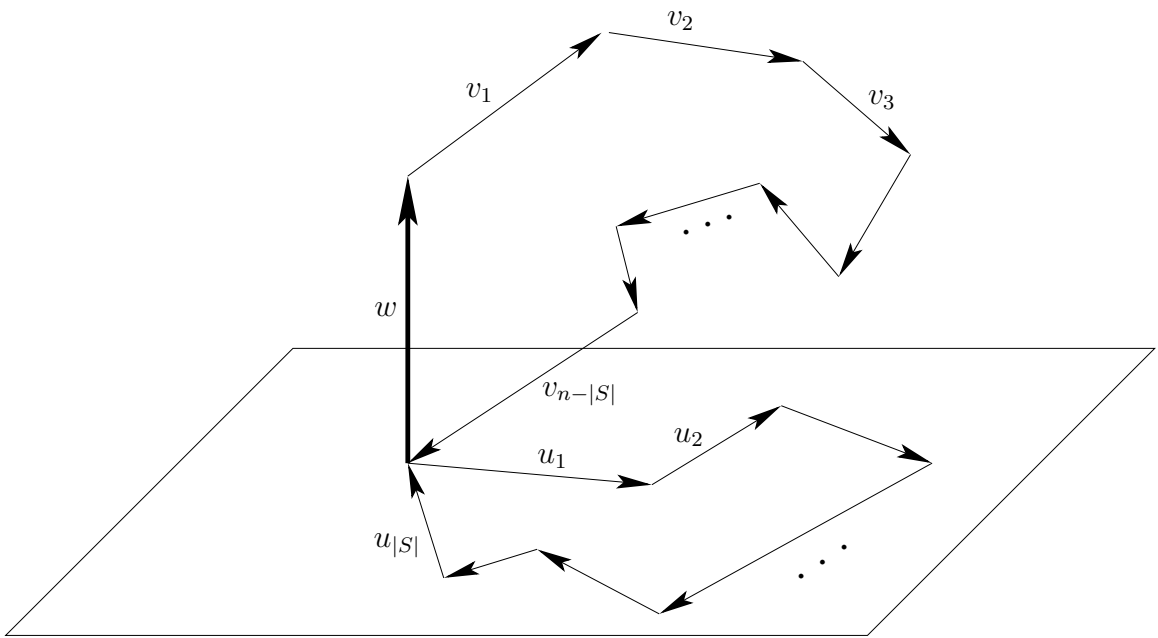}
\caption{
An element of $U_S$, represented by a spatial polygon with a distinguished
edge, and a planar polygon perpendicular to that edge.}
\label{fig:Uspoly}
\end{center}
\end{figure}
The first is the $n - |S| + 1$ sided polygon consisting of the vectors
$\{v_j\mid j\in S^c\}$ and $w$.  Each vector $v_j$ has length $\a_j$, and
$w$ has a variable length, always greater than or equal to $\sum_{i\in S}\a_i$.
This variable length is determined by the lengths of the edges in the second polygon, 
which consists of $|S|$ vectors $\{u_i\mid i\in S\}$, 
all contained in the plane perpendicular to $w$.
Note that this description also applies to the K\"ahler polygon space $M$ by taking $S=\emptyset$.

By setting $u_i = 0$ for all $i$ we get $M_S$,
the minimum of the Morse-Bott function $\Phi$ on $U_S$.
On the other hand, consider the submanifold of $U_S$ obtained by imposing
the extra condition that $\|w\| = \sum_{j\in S^c}\|v_j\|$.
Then the first of the two polygons is forced to be linear,
and we are left with $|S|$ vectors $\{u_i\}$ in the perpendicular plane
satisfying a certain norm condition
and adding to zero. 
Identifying this plane with $\C$ and dividing by the circle action rotating this plane, 
we obtain $\C P^{|S|-2}\cong X_S$,
the maximum of $\Phi$ on $U_S$.
Other critical points of $\Phi$ occur whenever the first polygon is linear,
which is possible for finitely many values of $\|w\|$.
\end{remark}

\begin{usproof}
Suppose given a point $[p,q]_{\R}\in U_S$,
and let
$$u_i = q_ip_i + p_i^*q_i^*\hs\hs\text{ for all }i\in S,$$
$$v_j = (q_j q_j^*)_0\hs\hs\text{ for all }j\in S^c,$$
$$w = \sum_{i\in S}(q_iq_i^*)_0-(p_i^*p_i)_0.$$
These vectors live in $i\cdot\sut\cong\sut^*\cong\R^3$, which is endowed with the
metric
$A\cdot B = \frac{1}{2}\operatorname{tr}AB$,
invariant under the coadjoint action.
With respect to this metric, we have the equalities
$\|(qq^*)_0\| = \half |q|^2$ and $\|(p^*p)_0\| = \half |p|^2$,
hence conditions (1), (2), and (4) are immediate consequences 
of the moment map equations.

To verify condition (3), note that the vectors $\{q_i\mid i\in S\}$
are all proportional over $\C$, which implies that the vectors $(q_iq_i^*)_0$
are positive scalar multiples of each other.
Furthermore, the moment map equation $p_iq_i=0$ implies that $(p_i^*p_i)_0$
is a non-positive scalar multiple of $(q_iq_i^*)_0$, therefore
$w = \sum(q_iq_i^*)_0-(p_i^*p_i)_0$ is proportional over $\R_+$ to
$(q_iq_i^*)_0$ for any $i\in S$.
Then $u_i\cdot w = \half\tr u_i w$ is a multiple of
$$\tr u_i(q_iq_i^*)_0 = \tr u_iq_iq_i^* = \tr p_i^*q_i^*q_iq_i^*
= |q_i|^2\tr p_i^*q_i^* = 0,$$
where the first equality comes from the fact that
$q_iq_i^*-(q_iq_i^*)_0$ is a scalar multiple of the identity,
and $\tr u_i = 0$.

To check condition (5), we first compute the norm of $u_i$:
\begin{eqnarray*}
\|u_i\|^2 &=& \half\tr u_i^2\\
&=& \half\tr(q_ip_ip_i^*q_i^* + p_i^*q_i^*q_ip_i)\\
&=& |q_i|^2|p_i|^2\\
&=& |q_i|^2 (|q_i|^2 - 2\a_i).
\end{eqnarray*}
Since all of the vectors
$\{(q_iq_i^*)_0 -(p_i^*p_i)_0\mid i\in S\}$ point in the same direction,
we have $$\|w\| = \sum_{i\in S}\|(q_iq_i^*)_0\| + \|(p_i^*p_i)_0\|
= \sum_{i\in S} \half|q_i|^2 + \half|p_i|^2
= \sum_{i\in S} |q_i|^2 - \a_i
= \sum_{i\in S} \sqrt{\a_i^2 + \|u_i\|^2}.$$

We have defined a map from $U_S$ to the moduli space of
vectors $\{u_i,v_j,w\}$ satisfying conditions (1)-(5),
and we claim that this map is a homeomorphism.
Since the source of this map is compact and the target is Hausdorff,
it is sufficient to show that the map is bijective.

Suppose given a collection of vectors $\{u_i,v_j,w\}\subs\sut$ 
satisfying conditions (1)-(5).  Using the adjoint action of $SU(2)$,
we may assume that $w$ is a positive scalar multiple of
$\footnotesize{\Big(\begin{array}{cc}
1&0\\
0& -1
\end{array}
\Big)}$.
By condition (3), this implies that for all $i\in S$, there exists $\la_i\in\C$
with $u_i =
\footnotesize{\Big(\begin{array}{cc}
0& \la_i\\
\bar{\la_i}& 0
\end{array}
\Big)}$. 
For $j\in S^c$, we choose $q_j\in\C^2$ with $(q_jq_j^*)_0=v_j$,
and observe that $q_j$ is unique up to the action of $U(1)^n$.
We know that for all $i\in S$, $(q_iq_i^*)_0$ must be a positive multiple of $w$,
hence there exist $a_i,b_i\in\C$ such that
$$q_i=\binom{a_i}{0}\hs\text{ and }\hs p_i = (0\hs b_i)$$
for all $i\in S$.
In order to have $u_i = q_ip_i + p_i^*q_i^*$ and $\half|q_i|^2-\half|p_i|^2=\a_i$, 
we must have $$a_i b_i = \la_i\hs\text{ and }\hs 
\half|a_i|^2-\half|b_i|^2=\a_i.$$
These equations uniquely define $a_i$ and $b_i$ up to the action of $U(1)^n$.
It follows from conditions (1)-(5) that 
$(p,q)\in\mur^{-1}(\a)\cap\muc^{-1}(0)$
and that $w=\sum_{i\in S}(q_iq_i^*)_0-(p_i^*p_i)_0$. 
This shows that our map is bijective, and thus completes the proof of Theorem \ref{us}.
\end{usproof}

\begin{remark}
Suppose that $S$ has only two elements; without loss of generality
we will assume that $S=\{1,2\}$.  Then forgetting $u_1$ and $u_2$
gives a diffeomorphism from $U_S$ to the ``vertical polygon space'' 
$VP(\a_3,\ldots,\a_n,\a_1+\a_2)$
defined in \cite{HK2}, shown to be diffeomorphic to a toric variety.
More generally with $S=\{1,\ldots,k\}$, given any two-element subset $T\subs S$, the subvariety
of $U_S$ given by the equations $u_i=0$ for all $i\in S\smallsetminus T$
is diffeomorphic to $VP(\a_{k+1},\ldots,\a_n,\sum_T\a_i)$.
\end{remark}
\end{section}

\begin{section}{Equivariant cohomology of $\mathbf{X}$}\label{s1}
We begin by defining representations
$$\rho_{ij}:K\to U(1)\hs\hs\hs\hs\text{ and }\hs\hs\hs\hs\rho_{\SO}:K\to SO(\sut)$$
by the formulae
$$\rho_{ij}\Ae = e_ie_j\hs\hs\hs\hs
\text{ and }\hs\hs\hs\hs\rho_{\SO}\Ae = \Ad(A).$$
Associated to these representations are the vector bundles
$$\L_{ij} = \Big(\mur^{-1}(\a)\cap\muc^{-1}(0)\Big)\times_{\rho_{ij}}\C\hs\hs\hs\hs
\text{ and }\hs\hs\hs\hs\mathcal{E}
=\Big(\mur^{-1}(\a)\cap\muc^{-1}(0)\Big)\times_{\rho_{\SO}}\sut.$$
Let $c_i = c_1(\L_{ii})\in H^2(X)$ and 
$p = p_1(\mathcal{E})\in H^4(X)$.

\begin{theorem}\label{Konno}{\em\cite{K2}}
The cohomology ring $H^*(X)$ is isomorphic to $\Q[c_1,\ldots,c_n,p]/\mathcal{I}$,
where $\mathcal{I}$ is generated by the following two families:
\begin{eqnarray*}&1)&\hs\hs p-c_i^2\hs\hs\text{ for all }\hs\hs i\in\on\\
&2)&\hs\hs\text{all elements of degree }2(n-2).
\end{eqnarray*}
\end{theorem}

The action of $S^1$ on the total spaces of $\L_{ij}$ and $\mathcal{E}$
given by the left action on $\mur^{-1}(\a)\cap\muc^{-1}(0)$ gives these bundles
an $S^1$-equivariant structure.  For the rest of this section,
we will use $c_i$ and $p$ to denote the {\em equivariant} characteristic classes
$$c_i = c_1(\L_{ii})\in H_{S^1}^2(X)\hs\hs\hs\hs\text{ and }\hs\hs\hs\hs
p = p_1(\mathcal{E})\in H_{S^1}^4(X).$$
Let $\K$ be the $S^1$-equivariant line bundle on $X$ obtained by pulling
back the weight $1$ line bundle over a point, and let $$x=c_1(\K)\in H_{S^1}^2(X).$$
We obtain the following result, extending Konno's work
to the equivariant context.

\begin{theorem}\label{main}
The equivariant cohomology ring $\hso(X)$ is isomorphic to $\Q[c_1,\ldots,c_n,p,x]/\mathcal{J}$,
where $\mathcal{J}$ is generated by the following two families:
\begin{eqnarray*}&1)&\hs\hs p-c_i^2\hs\hs\text{ for all }\hs\hs i\in\on\\
&2)&\hs\hs 
\Pd_{j\in\Lb}(c_j+\cns)\times\Pd_{i\in\Sb}(c_i+x)
\hs\hs\text{ for all }\hs\hs \emptyset\neq S\in\mathcal{S},
\end{eqnarray*}
where $m_S\in S$ and $n_S\in S^c$ are the minimal elements of the two sets,
$\Sb = S\smallsetminus\{m_S\}$, and $\Lb=S^c\smallsetminus\{n_S\}$.
\end{theorem}

\begin{remark}  Konno observes that the quotient map from the abstract
polynomial ring $\Q[c_1,\ldots,c_n,p]$ to $H^*(X)$ is precisely the Kirwan
map $$\k:\hg(\cot)\to H^*(X)$$ induced by the inclusion
$\mur^{-1}(\a)\cap\muc^{-1}(0)\hookto\cot$.  
Theorem \ref{Konno} can be interpreted as saying
that the Kirwan map for hyperpolygon spaces is surjective, with kernel $\mathcal{I}$.
Likewise, Theorem \ref{main} asserts that the $S^1$-equivariant Kirwan map
$$\kso:\hgso(\cot)\to\hso(X)$$ is surjective, with kernel $\mathcal{J}$.
The analogous map for K\"ahler reductions
is known to always be surjective \cite{Ki}, but in the hyperk\"ahler case
the problem remains open.
\end{remark}

\begin{remark}\label{meaning}
The second family of relations in Theorem~\ref{main} has a geometric
interpretation. 
In the course of the proof of Theorem~\ref{main} it will be shown that
the class \(-\half(c_j + c_{n_s})\) is represented by the divisor given
by the points \([p,q] \in X\) where $q_j$ and $q_{n_s}$ are
straight. Hence the product \(\prod_{i \in \Lb}(c_j + c_{n_s})\) is
supported on the subvariety of points where $S^c$ is
straight. Similarly, it is shown in the proof that the class \(c_i +
x\) is represented by the divisor given by the condition \(p_i = 0,\)
hence the class
\[
\Pd_{j\in\Lb}(c_j+\cns)\times\Pd_{i\in\Sb}(c_i+x)
\]
is supported on the subvariety of points $[p,q]$ where a long subset
$S^c$ is straight, and {\em all but one} $p_i$ for \(i \in S\) is
zero. 
We will show, using Theorem \ref{fixed}, 
that this subvariety is disjoint from the fixed point set $X^{S^1}$,
and that this implies that its cohomology class is trivial.
\end{remark}

Before proceeding with the proof of Theorem \ref{main}, we collect some
preliminary results regarding the relationship between ordinary and equivariant cohomology.
The first result that we need is known as {\em equivariant formality},
proven for compact manifolds in \cite{Ki}, which we adapt to our situation
in Proposition \ref{formality}.

\begin{proposition}\label{formality}
Let $X$ be a symplectic
manifold, possibly noncompact but of finite topological type.
Suppose that $X$ admits a hamiltonian circle action,
and that the moment map
is proper and bounded below.  Then
$H^*_{S^1}(X)$
is a free module over $H^*_{S^1}(pt)$.
\end{proposition}

\begin{proof}
Because $\Phi$ is a moment map, it is a Morse-Bott function
such that all of the critical submanifolds and their normal bundles
carry almost complex structures.  Thus we get a Morse-Bott decomposition
of $X$ into even-dimensional $S^1$-invariant submanifolds.
This tells us, as in \cite{Ki},
that the spectral sequence associated to the fibration
$X\hookto EG\times_G X\to BG$ collapses, and we get the desired result.
\end{proof}

Consider the following commuting square of maps, where
$\phi$ and $\psi$ are each given by setting $x$ to zero.

$$
\xymatrix{
0 \ar[r] & \mathcal{J} \ar[r]\ar[d]_{\phi} & \hgso(\cot) \ar[r] \ar[d]_{\phi}
& \hso(X) \ar[d]_{\psi}\\
0 \ar[r] & \mathcal{I} \ar[r] & \hg(\cot) \ar[r] & H^{*}(X) }
$$

\vspace{\baselineskip}
\noindent Our moment map $\Phi:X\to\R$ is proper and bounded below \cite[1.3]{HP}, therefore
Proposition \ref{formality} has the following consequences.

\begin{corollary}\label{surjective}
The equivariant Kirwan map $\kso$ is surjective.
\end{corollary}

\begin{proof}
Suppose that $\gamma\in\hso(X)$ is a homogeneous class of minimal degree
that is {\em not}
in the image of $\kso$, and choose a class
$\eta\in\phi^{-1}\k^{-1}\psi(\gamma)$.
Then by Proposition \ref{formality}, $\kso(\eta)-\gamma = x\delta$ for some $\delta\in\hso(M)$,
and therefore $\delta$ is a class of lower degree that is not in the image of $\kso$.
\end{proof}

\vspace{-\baselineskip}
\begin{corollary}\label{enough}
Suppose that $\mathcal{J}\subs\ker\kso$ and $\phi(\mathcal{J}) = \mathcal{I}$.
Then $\mathcal{J}=\ker\kso$.
\end{corollary}

\begin{proof}
Suppose that $a\in\ker\kso\smallsetminus\mathcal{J}$ is a homogeneous class of minimal degree,
and choose $b\in\mathcal{J}$ such that $\phi(a-b)=0$.  Then $a-b = cx$ for some
$c\in\hgso(\cot)$.  By Proposition \ref{formality},
$cx\in\ker\kso\impl c\in\ker\kso$, hence $c\in\ker\kso\smallsetminus\mathcal{J}$ 
is a class of lower degree than $a$.
\end{proof}

\vspace{-\baselineskip}
\begin{corollary}\label{injectivity}
Let $E$ be an $S^1$-equivariant vector bundle on $X$, with an equivariant section
$s$ such that the zero set of $s$ is disjoint from the fixed point set $X^{S^1}$.
Then $e(E)=0\in\hso(X)$.
\end{corollary}

\begin{proof}
Consider the bundle $\K^{\oplus 2(n-3)}$ restricted to $X\smallsetminus X^{S^1}$.
An equivariant section of this bundle is equivalent to an ordinary section of the
induced bundle over the quotient $\left(X\smallsetminus X^{S^1}\right)\big/S^1$,
and by degree considerations we can find a nonvanishing section.
Hence $\K^{\oplus 2(n-3)}$ has an equivariant section $t$ over $X$ with zero set
supported on $X^{S^1}$.  Then $s\oplus t$ is a nonvanishing section of
$E\oplus \K^{\oplus 2(n-3)}$, hence
$$0 = e\left(E\oplus \K^{\oplus 2(n-3)}\right) 
= e(E)\cdot e\left(\K^{\oplus 2(n-3)}\right)
=x^{2(n-3)}e(E).$$
Then by Proposition \ref{formality}, $e(E)=0$.
\end{proof}

\vspace{-\baselineskip}
\begin{remark}
Corollary \ref{injectivity} is a weak form of the statement that the restriction
map from $\hso(X)$ to $\hso(X^{S^1})$ is injective, proven for compact $X$ in \cite{Ki}.
\end{remark}

\begin{mainproof}
Corollary \ref{surjective} tells us that the characteristic classes
$c_1,\ldots,c_n,p,x\in\hso(X)$ generate the ring, and Corollary \ref{enough}
tells us that it is enough to prove two statements:  the first is that 
$\mathcal{J}\subs\ker\kso$, i.e. that the elements of $\mathcal{J}$ are
indeed relations in $\hso(X)$, and the second is that 
$\phi(\mathcal{J})=\mathcal{I}$.
We begin by proving that $\mathcal{J}\subs\ker\kso$, following the approach
outlined in Remark \ref{meaning}.

In the nonequivariant context, 
Konno shows that $\mathcal{E}\cong\L_{ii}\oplus\R$
as a real vector bundle for all $i$ \cite{K2}.  This implies that
$$p = p_1(\mathcal{E}) = -c_2(\mathcal{E}\otimes\C)
= -c_2(\L_{ii}\oplus\L_{ii}^*) = c_1(\L_{ii})^2 = c_i^2\in H^*(X),$$ 
and therefore $p-c_i^2\in\ker\k$.
This argument is adaptable to the equivariant context without any modifications,
hence we will not include it here.  

Consider the function $$\tilde{s}_{ij}:\mur^{-1}(\a)\cap\muc^{-1}(0)\to\C$$
given by $\tilde{s}_{ij}(p,q) = \operatorname{det}(q_i q_j)$,
where $(q_i q_j)$ is considered to be a $2\times 2$ matrix \cite{K2}.
This function is $S^1$-invariant, and $K$-equivariant with respect to $\rho_{ij}$,
and therefore defines an $S^1$-equivariant 
section $s_{ij}$ of $\L_{ij}^*$ (the dualization is a consequence
of the fact that the action of $K$ on $\mur^{-1}(\a)\cap\muc^{-1}(0)$ is a right action).
The vanishing set of $s_{ij}$ is the divisor 
$$Z_{ij} = \big\{[p,q]\in X\bigmid q_i\text{  is proportional to  }q_j\big\},$$
hence $c_1(\L_{ij}^*) = -\half(c_i+c_j)$
is represented in equivariant Borel-Moore homology by the divisor $Z_{ij}$.
Now consider the function $$\tilde{t}_i:\mur^{-1}(\a)\cap\muc^{-1}(0)\to\C$$
given by 
$$\tilde{t}_i(p,q) = 
\begin{cases} 
p_i^{(2)}/q_i^{(1)} & \text{if $q_i^{(1)}\neq 0$}\\
-p_i^{(1)}/q_i^{(2)} & \text{if $q_i^{(2)}\neq 0$,}
\end{cases}$$
where $$q_i = \left(\begin{array}{c}q_i^{(1)}\\q_i^{(2)}\end{array}\right)
\hs\text{ and }\hs p_i = (p_i^{(1)} p_i^{(2)}).$$
This function is well-defined by the fact that $\mc(p,q)=0$. 
It is $S^1$-equivariant 
with respect to the weight 1 action of $S^1$ on $\C$,
and it is $K$-equivariant with respect to $\bar\rho_{ii}$.
Thus it descends to an equivariant section $t_i$ of $\L_{ii}\otimes \K$,
vanishing on the divisor
$$W_i = \big\{[p,q]\in X\bigmid p_i=0\big\},$$
so $c_i+x=c_1(\L_{ii}\otimes \K)$ is represented by the divisor $W_i$.

Consider a nonempty short subset $S\in\mathcal{S}$, and
define the vector bundle
$$E_S = \bigoplus_{j\in\Lb}\L_{jn_S}^*\hs\hs\oplus\hs\hs
\bigoplus_{i\in\Sb}\L_{ii}\otimes \K$$
with equivariant Euler class
$$\left(-1/2\right)^{|\Lb|}
\Pd_{j\in\Lb}(c_j+\cns)\times\Pd_{i\in\Sb}(c_i+x).$$
Then $\left(\oplus_{i\in\Sb}t_i\right)\oplus_{j\in\Lb}s_{jn_S}$
is a section of $E_S$ vanishing only on the cycle 
$$Z_S = \bigcap_{i\in\Sb}Z_i\hs\hs\cap\hs\hs\bigcap_{j\in\Lb}Z_{jn_S}$$
consisting of points $[p,q]$
such that $q_j$ is proportional to $q_{n_S}$ for all $j\in S^c$,
and $p_i=0$ for all $i\in\Sb$.
We would like to show that $e(E_S)=0\in\hso(X)$, and
by Corollary \ref{injectivity}, it will suffice to show that $Z_S$ is disjoint from
$$X^{S^1} = M\cup\bigcup_{T\in\mathcal{S}}X_T.$$
Since $S^c$ is a long subset that is straight in $Z_S$, we have $Z_S\cap M=\emptyset$
by Theorem \ref{stability}.
We must now show that $Z_S$ also does not intersect any $X_T \subset
X^{S^1}$. We begin with the observation that for each \(X_T, T \in
\mathcal{S}'\), 
\begin{equation}\label{nonzerop}
[p,q] \in X_T \Rightarrow \text{ at least {\em two} of the vectors in
} \{p_i \mid i\in T\} \text{ are {\em nonzero}.} 
\end{equation}
This follows from the description of $X_T$ given in
Theorem~\ref{fixed} and the complex moment map conditions. Now let \(T
\in {\mathcal S}'\) be a short subset.
If \(T \nsubseteq S,\) then by the descriptions of $X_T$
and $Z_S$, we may conclude that any point \([p,q] \in Z_S \cap X_T\)
must have the long subset
\(T \cup S^c\) straight, and \(p_j=0\) for all \(j \in T^c.\) This
means that $(p,q)$ is unstable by Theorem~\ref{stability}, so \(Z_S \cap
X_T = \emptyset.\) On the other hand, if \(T \subseteq S,\) then we
may similarly conclude that for any point \([p,q] \in Z_S \cap X_T,\) 
the long subset $T^c$ is straight, and {\em at most one} element in
\(\{p_i \mid i \in T\}\) is nonzero. This contradicts the
observation~\eqref{nonzerop} above, so \(Z_S \cap X_T = \emptyset.\)

We have now proven that $\mathcal{J}\subs\ker\kso$, and it remains
only to show that $\phi(\mathcal{J})=\mathcal{I}$.
In other words (by Theorem \ref{Konno}), we must show that the set
$$\left\{\Pd_{j\in\Lb}(c_j+\cns)\times\Pd_{i\in\Sb}c_i\hs\hs\Bigg |\hs\hs
\emptyset\neq S\in\mathcal{S}\right\}$$
spans all monomials of degree $2(n-2)$ in the ring
$\Q[c_1,\ldots,c_n]\big/\big<c_i^2-c_1^2\hs\big{|}\hs i\in\twn\big>$.

Let $\bk = \half(c_1+\ck)$ for all $k$, so that $c_k = 2\bk-b_1$.
The relations $c_k^2 = c_1^2$ translate to $b_k^2 = b_1b_k$ for all $k$.
Let $$v_S = \frac{(-1)^n}{2^{|\Lb|}}\Pd_{j\in\Lb}(c_j+\cns)\times\Pd_{i\in\Sb}c_i
= (-1)^n\pjl (\bj+\bns-b_1)\times\pis(2\bi-b_1),$$
and let
$$\ba = (-1)^{|A|} b_1^{n-2-|A|} \prod_{k\in A} b_k$$ for all $A \subsetneq \twn$.
Then $\{b_A\}$
is a basis for
the $(n-2)^{\text{nd}}$ graded piece of the ring
$$\Q[b_1,\ldots,b_n]\big/\big<b_k^2-b_1b_k\hs\big{|}\hs k\in\twn\big>,$$
hence we need to show that each element $b_A$ can be expressed as a linear combination
of the elements $\{v_S\mid S\in\sa\}$.

\begin{claim}\label{vs}
We have the following expression for $v_S$ in terms of the basis $\{b_A\}$:
$$\vs =
\begin{cases}
\displaystyle{\sum_{\ss{\Lb\hs\subseteq A \\ m_S\notin A}}} \hs 2^{|\asb|} \hs\ba
& \text{if}\hspace{10pt}1 \in S^c;\vspace{.5cm}\\
\displaystyle{\sum_{S^c \nsubseteq A}} \hs 2^{|\asb|} 
\hs\ba & \text{if}\hspace{10pt}1 \in S.\\
\end{cases}
$$
\end{claim}

\begin{proof}
Any degree $n-2$ monomial
in $b_1,\ldots,b_n$ is equal to $(-1)^{|A|}b_A$, where $A$ is the set of $k>1$ such that
$b_k$ appears in the monomial.
Expanding $\vs$, we need to count (with sign) the occurrence of $b_A$ for each $A$.
In most cases we find that there is no cancellation, and the claim is straightforward.
The most difficult case occurs when
$1\in S$ (therefore $n_S=1$) and $\ns\in A$; 
in this case the number of times (with multiplicity) that $b_A$
occurs in $\vs$ is 
\begin{eqnarray*}
&& (-1)^n(-1)^{|A|}(-1)^{|\acsb|}\hs 2^{|\asb|}\sum_{E\subsetneq\Ac\cap\Lb}(-1)^{|E|}
\vspace{.3cm}\\
\vspace{.3cm}&=& (-1)^n(-1)^{|A|}(-1)^{|\acsb|}\hs 2^{|\asb|}
\left((1-1)^{|\Ac\cap\Lb|}-(-1)^{|\Ac\cap\Lb|}\right)\vspace{.3cm}\\
\vspace{.3cm}&=& (-1)^{n+|A|+|\acsb|+|\Ac\cap\Lb|+1}\hs 2^{|\asb|}\\
\vspace{.3cm}&=& (-1)^{2n}\hs 2^{|\asb|}\\
\vspace{.3cm}&=& \hs 2^{|\asb|}.
\end{eqnarray*}
(When we write $A^c$, we mean the complement of $A$ inside of $\{2,\ldots,n\}$.)
We leave the remaining cases to be checked by the reader.
\end{proof}

\vspace{-\baselineskip}
\begin{claim}\label{ws}
Suppose that $1\in S$.
Let $S_0 = S$, and for $1\leq k\leq |S|$, let $S_k = S_{k-1}\setminus\{m_{S_{k-1}}\}$
(i.e. $S_k$ consists of the $|S|-k$ largest elements of $S$).
Then $v_S + \displaystyle{\sum_{k=1}^{|S|-1}} \hs 2^{k-1}\hs v_{S_k} 
= \displaystyle{\sum_A}\hs 2^{|\asb|}\hs b_A$.
\end{claim}

\begin{proof}
We proceed by induction to show that
$$v_S + \displaystyle{\sum_{k=1}^l} \hs 2^{k-1}\hs v_{S_k} =
\displaystyle{\sum_A} \hs 2^{|\asb|}\hs b_A - 
\hs 2^l\hs\displaystyle{\sum_{\overline{S^c_{l+1}}\subs A}}
2^{|\asb_l|}\hs b_A;$$
the case $l=|S|-1$ is the statement of the claim.
The base case $l=0$ follows from Claim \ref{vs},
together with the observation that $\overline{S^c_{1}}=S^c$.
More generally, for all $l\geq 1$, we have $\overline{S^c_{l+1}} = S^c\cup\{m_{S_1},\ldots,m_{S_l}\}$.
Then 
\begin{eqnarray*}
v_S + \displaystyle{\sum_{k=1}^{l+1}} \hs 2^{k-1}\hs v_{S_k} &=&
v_S + \displaystyle{\sum_{k=1}^l}\hs 2^{k-1}\hs v_{S_k} + 2^l\hs v_{S_{l+1}}\\
&=& \displaystyle{\sum_A} \hs 2^{|\asb|}\hs b_A - 2^l\hs\displaystyle{\sum_{\overline{S^c_{l+1}}\subs A}}
2^{|\asb_l|}\hs b_A + 2^l\hs \sum_{\ss{\overline{S^c_{l+1}}\subs A \\ m_{S_{l+1}}\notin A}}
2^{|\asb_{l+1}|}\hs b_A
\end{eqnarray*}
by the inductive hypothesis and Claim \ref{vs}.
Using the fact that
$\asb_{l+1}=\asb_l$ when $m_{S_{l+1}}\notin A$,
this is equal to
$$\sum_A\hs 2^{|\asb|}\hs b_A - 2^l\hs\sum_{\overline{S^c_{l+1}}\cup\{m_{S_{l+1}}\}\subs A}
2^{|\asb_l|}.$$
Finally, since $|\asb_{l+1}|=|\asb_l|-1$ when $m_{S_{l+1}}\in A$,
this reduces to
$$\sum_A \hs 2^{|\asb|}\hs b_A - 2^{l+1}\hs\displaystyle{\sum_{\overline{S^c_{l+2}}
\subs A}}
2^{|\asb_{l+1}|}\hs b_A,$$ thus proving our claim.
\end{proof}

For all short subsets $T$ containing $1$,
let $w_T=\displaystyle{\sum_A}2^{|A\cap \bar{T}|}\hs b_A$, 
which by Claim \ref{ws} is expressible as a linear combination of elements
of the set $\{v_S\mid \emptyset\neq S\in\mathcal{S}\}$.
Let $$x_S = 
\begin{cases}
\displaystyle{\sum_{1\in T\subs S}}(-1)^{|S|+|T|}w_T & \text{if $1\in S$,}\\
v_S & \text{if $1\in S^c$.}
\end{cases}$$
Our last task will be to prove that the transition matrix $Q$ 
taking the basis $\{b_A\}$ to 
the set $\{x_S\}$
is upper triangular with ones on the diagonal, and therefore invertible.
In order to make sense of ``the diagonal," 
we must first give an explicit bijection between the set of proper subsets of $\twn$
and the set of nonempty short subsets of $\on$.
We do this as follows:  given $A\subsetneq\twn$, let 
$$S(A) = \begin{cases}
\Ac &\text{ if $\Ac$ is short,}\\
\on\setminus\Ac = A\cup\{1\} &\text{ if $\Ac$ is long.}
\end{cases}$$
The rows of $Q$ will be indexed by $A$, 
and the sets will appear in lexicographic order within cardinality
class.  For example, when $n=4$, the order of the rows will be 
$\emptyset$, $\{2\}$, $\{3\}$, $\{4\}$,
$\{2,3\}$, $\{2,4\}$, $\{3,4\}$.  
The columns will be indexed by $S$ according to the bijection described above.

\begin{claim}\label{tri}
The matrix $Q$ is lower triangular with ones on the diagonal.
\end{claim}

\begin{proof}
First consider a column corresponding to a short subset $S$ that does {\it not} contain 1.
The entries in this column correspond to the coefficient of $b_A$ in $x_S = v_S$.
From Claim \ref{vs}, we see that $b_A$ appears in $v_S$ only if $\Lb\subs A\subs \Lb\cup\Sb$,
and if so it appears with a coefficient of $2^{|\asb|}$.
Since $1\notin S$, we have $\Lb = S^c\setminus\{1\} = \twn\setminus S$.
The diagonal entry corresponds to the set $A = \twn\setminus S = \Lb$, therefore
in this row we get the number $2^{|\asb|} = 2^{|\Lb\cap\Sb|} = 1$.
Since the set $A$ corresponding to a given row can never contain the set $B$ corresponding
to a 
lower row, the rows above the diagonal fail to satisfy the condition $\Lb\subs A$,
and we get all zeros.

Now consider a column corresponding to a short subset $S$ that {\it does} contain $1$.
In this case, the coefficient of $b_A$ in $x_S$ is
$$(-1)^{|S|}\displaystyle{\sum_{1\in T\subs S}}(-1)^{|T|}2^{|A\cap\bar{T}|}.$$
The diagonal entry corresponds to the set $A=\Sb$, and we get
\begin{eqnarray*}
(-1)^{|S|}\displaystyle{\sum_{1\in T\subs S}}(-1)^{|T|}2^{|\bar{T}|}
&=& (-1)^{|\Sb|}\displaystyle{\sum_{1\in T\subs S}}(-2)^{|\bar{T}|}\\
&=& (-1)^{|S|}(1-2)^{|\Sb|} = 1.
\end{eqnarray*}
Any row above the diagonal corresponds to a set $A$ which does not contain $\Sb$.
Choose an element $l\in\Sb\setminus A$.  Then
\begin{eqnarray*}
(-1)^{|S|}\displaystyle{\sum_{1\in T\subs S}}(-1)^{|T|}2^{|A\cap\bar{T}|}
&=& (-1)^{|S|}\displaystyle{\sum_{l\in T}}(-1)^{|T|}2^{|A\cap\bar{T}|}
+ (-1)^{|S|}\displaystyle{\sum_{l\notin T}}(-1)^{|T|}2^{|A\cap\bar{T}|}\\
&=& (-1)^{|S|}\displaystyle{\sum_{l\notin T}}\left[(-1)^{|T|}2^{|A\cap\bar{T}|}
+ (-1)^{|T\cup\{l\}|}2^{|A\cap\bar{T}|}\right]\\
&=& 0.
\end{eqnarray*}
\end{proof}

Claim \ref{tri} tells us that each $b_A$ can be expressed as a linear
combination of elements of the form $x_S$, and therefore of elements
of the form $v_S$.  This lets us conclude that $\phi(\mathcal{J})=\mathcal{I}$,
and thereby completes the proof of Theorem \ref{main}.
\end{mainproof}
\end{section}

\begin{section}{The cohomology ring of a core component}\label{last}

In this section we compute the $S^1$-equivariant 
and ordinary cohomology rings of the core component
$U_S$ corresponding to a short subset $S\subs\on$.  Since $U_S$ is the closure of a cell
in an even-dimensional equivariant cellular decomposition of $X$, the restriction
map $\hso(X)\to\hso(U_S)$ is surjective.  In particular, $\hso(U_S)$ 
is generated by restrictions
of the Kirwan classes $c_1,\ldots,c_n,x$.  
For our presentation, it will be convenient
to assume that $1\in S$, and to work with the classes 
$b_k = \half(c_1+c_k)$ introduced in Section
\ref{s1}.  
With respect to these generators, we obtain the following result.

\begin{theorem}\label{eqcore}
The equivariant cohomology ring $\hso(U_S)$ is isomorphic to
$\Q[b_1,\ldots,b_n,x]/\mathcal{J}_S,$
where $\mathcal{J}_S$ is generated by the following four families:
\begin{eqnarray*}&1)&\hs\hs b_1 - b_i \hs\hs\text{ for all }\hs\hs i\in S\\
&2)&\hs\hs
b_j(b_1 - b_j)
\hs\hs\text{ for all }\hs\hs j\in S^c\\
&3)&\hs\hs \prod_{j\in R}b_j\hs\hs\text{ for all }\hs\hs R\subs S^c\text{ such that }R\cup S
\text{ is long}\\
&4)&\hs\hs (b_1 + x)^{|S|-1}\cdot\frac{1}{b_1}\left(\prod_{j\in L}(b_j-b_1)-\prod_{j\in L}b_j\right)
\hs\hs\text{ for all long subsets }L\subs S^c.
\end{eqnarray*}
\end{theorem}

\begin{corollary}\label{ordcore}
The ordinary cohomology ring $H^*(U_S)$ is isomorphic to
$\Q[b_1,\ldots,b_n]/\mathcal{I}_S,$
where $\mathcal{I}_S$ is generated by the following four families:
\begin{eqnarray*}&1)&\hs\hs b_1 - b_i \hs\hs\text{ for all }\hs\hs i\in S\\
&2)&\hs\hs
b_j(b_1 - b_j)
\hs\hs\text{ for all }\hs\hs j\in S^c\\
&3)&\hs\hs \prod_{j\in R}b_j\hs\hs\text{ for all }\hs\hs R\subs S^c\text{ such that }R\cup S
\text{ is long}\\
&4)&\hs\hs b_1^{|S|-2}\prod_{j\in L}(b_j-b_1)
\hs\hs\text{ for all long subsets }L\subs S^c.
\end{eqnarray*}
\end{corollary}

\begin{remark}\label{integer}
Each of these relations has a geometric interpretation.
For $i\in\{1,\ldots,n\}$, it is possible to construct a line bundle on $X$ with equivariant
Euler class $b_i - b_1$
which has a section supported on the locus
where $q_1q_1^*$ and $q_iq_i^*\in\R^3$ point in opposite directions.  
Since this locus is disjoint
from $U_S$ when $i\in S$, we have $$1)\hs\hs b_i = b_1 \in \hso(\us)\hs\text{ for all }\hs i\in S.$$
Similarly, we showed in the proof of Theorem \ref{main}
that $-b_j = -\half(c_1+c_j)$ is represented by the divisor $Z_{1j}$
on which $q_1q_1^*$ and $q_iq_i^*\in\R^3$ point in the same direction.
Then by the previous reasoning, we obtain 
$$2)\hs\hs b_j(b_1 - b_j) = 0\in\hso(\us)\hs\text{ for all }\hs j\in S^c.$$
Recall from Section \ref{s1} that for any $R\subs S^c$,
the cohomology class $(-1)^{|R|}\prod_{j\in R}b_j$ is represented
by the subvariety $Z_R\subs X$ of points with $q_j$ proportional to $q_1$
for all $j\in R$.
When restricted to $U_S$, this becomes $U_S\cap U_{R\cup S}$,
the unstable manifold for the critical locus $X_{R\cup S}\cap U_S$
of the Morse-Bott function $\Phi|_{U_S}$.
In particular, we have $$3)\hs\hs \prod_{j\in R}b_j=0\in\hso(\us)\text{ if }R\cup S\hs\text{ is long}.$$
To understand the fourth family of relations, recall from Section \ref{s1} that
$$b_1+x = 2b_i-b_1+x = c_i+x \in \hso(U_S)$$ 
is represented by the divisor $W_i$ of points
with $p_i = 0$ for any $i\in S$.  In particular, the class 
$(b_1+x)^{|S|-1}$ 
is represented by the subvariety of points in $U_S$ on which $p_i = 0$
for all $i\in\Sb$, which is equal to $M_S$ by the complex moment map condition.
Hence the fourth family of generators of $\mathcal{J}_S$ (or of $\mathcal{I}_S$)
can be interpreted geometrically as $(b_1+x)^{|S|-1}$ (respectively $b_1^{|S|-1}$
in the nonequivariant case) times
classes that vanish in $\hso(M_S)$ (see Lemma \ref{pols}).
%
\end{remark}

\begin{eqcoreproof}
Let $\phi:\Q[b_1,\ldots,b_n,x]\to\hso(\us)$ denote the composition of the Kirwan
map with restriction to $\us$.  Our claim is that $\ker\phi = \mathcal{J}_S$.
For every short subset $T$ containing $S$,
let $$\phi_T:\Q[b_1,\ldots,b_n,x]\to\hso(\mtus)$$ denote
the composition of the Kirwan map with restriction to $\mtus$,
and let $$J_T = \ker\phi_T.$$
Similarly, let $$\phi_{\emptyset}:\Q[b_1,\ldots,b_n,x]\to\hso(M_S)$$
be the natural map, and let
$$J_{\emptyset} = \ker\phi_{\emptyset}.$$
The kernel of the restriction map $\hso(U_S)\to\hso(U_S^{S^1})$
to the fixed point set of $U_S$ is a torsion module over $\hso(pt)$ \cite[3.5]{AB}, 
and Proposition \ref{formality}
tells us that $\hso(U_S)$ is a free $\hso(pt)$-module. 
Hence the restriction map is injective,
and we have $$\ker\phi = \ker\phi_{\emptyset}\cap\bigcap_{\ts}\ker\phi_T.$$
We know that $\mtus\cong\C P^{|S|-2}$ for all short $\ts$,
therefore $$\hso(\mtus)\cong\Q[h,x]/h^{\sz}.$$
Furthermore, we know that for all $i\in T$, the restriction of $b_i + x$ 
to $\hso(U_T)$ is represented by the divisor $W_i\cap U_T$ (see Remark \ref{integer}),
and therefore restricts to the class of a hyperplane on $\mtus$.
Hence $\phi_T(b_i+x)=h$ for all $i\in T$.
On the other hand, for $j\in T^c$, the class $b_j$ 
is represented by the divisor $Z_{1j}$ on $X$, which is disjoint from $\mtus$,
hence $\phi(b_j) = 0$ for all $j\in T^c$.
Thus we conclude that
$$\ker\phi_T = \<b_1-b_i, \hs b_j, (b_1+x)^{\sz}\mid i\in T, j\in T^c\>.$$

\begin{lemma}\label{jt}
The intersection 
$\displaystyle{
\bigcap_{\ts}}\ker\phi_T$
is equal to
$$\<b_1-b_i, \hs b_j(b_1 - b_j), 
\displaystyle{\prod_{j\in R}b_j},
(b_1+x)^{\sz}
\hs\hs\bigg |\hs\hs i\in S, j\in S^c, R\cup S\text{ long}\>.$$
\end{lemma}

\begin{proof}
First, since the variable $x$ appears only in the generator $(b_1+x)^{\sz}$,
which is contained in every ideal on both sides of the equation, we may reduce
the problem to showing that
\begin{equation}\label{reduce}
\bigcap_{\ts}\big\langle b_1-b_i, \hs b_j\bigmid i\in T, j\in T^c\big\rangle = 
\<b_1-b_i, \hs b_j(b_1-b_j), \prod_{j\in R}b_j
\hs\hs\bigg |\hs\hs i\in S, j\in S^c, R\cup S\text{ long}\>
\end{equation}
in the ring $\Q[b_1,\ldots,b_n]$.
Both ideals cut out the (reducible) variety 
$$\bigcup_{\ts}Y_T\subs\operatorname{Spec}\Q[b_1,\ldots,b_n],$$
where $$Y_T = \big\{(z_1,\ldots,z_n\bigmid
z_i=z_1\hs\forall i\in S, z_j=0\hs\forall j\in S^c\big\}.$$
The left hand side of Equation \eqref{reduce} is an intersection of prime ideals,
and is therefore radical.
Thus by Hilbert's Nullstellensatz, it is sufficient to prove that
the right hand side of Equation \eqref{reduce} is radical.
This involves showing that the ideal is saturated, with
Hilbert polynomial equal to the constant $\#\{\text{short }\ts\}$.

The degree $k$ piece of the quotient
$$\Q[b_1,\ldots,b_n]/\<b_1-b_i, \hs b_j(b_1-b_j)
\mid i\in S, j\in S^c\>$$
has a basis of elements of the form
$$b_1^{e_1}\prod_{j\in S^c}b_j^{e_j},$$
where $e_j \in\{0,1\}$ for all $j>0$, and $e_1+\sum_{j\in S^c}e_j = k$.
The subset of these elements with the property that
$S\cup\{j\mid e_j = 1\}$ is short descends to a basis for the
degree $k$ part of the ring
$$\Q[b_1,\ldots,b_n]\bigg/\<b_1-b_i, \hs b_j(b_1-b_j), \prod_{j\in R}b_j
\hs\hs\bigg |\hs\hs i\in S, j\in S^c, R\cup S\text{ long}\>,$$
hence our ideal has the desired Hilbert polynomial.
It is also clear from this description that if an element
$a$ of the quotient ring is nonzero, so is $b_1^d\cdot a$ for any $d\geq 0$,
hence our ideal is saturated.
\end{proof}

It now remains to show that
$$\mathcal{J}_S = \<b_1-b_i, \hs b_j(b_1 - b_j), \prod_{j\in R}b_j,
\hs (b_1+x)^{\sz} \hs\hs\bigg |\hs\hs i\in S, j\in S^c, R\cup S\text{ long}\> \cap \ker\phi_{\emptyset}.$$
The fact that $\mathcal{J}_S$ is contained in the intersection is clear.
To show the opposite containment, consider an element
$$a + \eta\cdot (b_1+x)^{\sz}\in \<b_1-b_i, \hs b_j(b_1 - b_j), \prod_{j\in R}b_j,
\hs (b_1+x)^{\sz} \hs\hs\bigg |\hs\hs i\in S, j\in S^c, R\cup S\text{ long}\>,$$
with $$a\in \<b_1-b_i, \hs b_j(b_1 - b_j), \prod_{j\in R}b_j\hs\hs\bigg |\hs\hs
i\in S, j\in S^c, R\cup S\text{ long}\>,$$
and suppose that we also have
$$a + \eta\cdot (b_1+x)^{\sz}\in\ker\phi_{\emptyset}.$$

\begin{lemma}\label{pols}{\em \cite{HK2}}
The kernel of $\phi_{\emptyset}$
is equal to 
$$\<b_1-b_i, \hs b_j(b_1-b_j), \prod_{j\in R}b_j, 
\hs (b_1+x)^{\sz}b_1^{-1}\left(\prod_{j\in L}(b_j-b_1)-\prod_{j\in L}b_j\right)\>,$$
where $i\in S, j\in S^c$, and $R,L\subs S^c$, with $R\cup S$ and $L$ both long.
\end{lemma}

Lemma \ref{pols} tells us that $a\in \ker\phi_{\emptyset}$,
therefore $$\eta\cdot (b_1+x)^{\sz}\in\ker\phi_{\emptyset}.$$
But $(b_1+x)^{\sz}$ is represented in $\hso(U_S)$ by the subvariety
$M_S$ (see Remark \ref{integer}), hence 
$$0 = \phi_{\emptyset}(\eta\cdot (b_1+x)^{\sz})
= \phi_{\emptyset}(\eta)\cdot e(M_S),$$
where $e(M_S)$ is the equivariant Euler class of the normal bundle
to $M_S$ inside of $U_S$.  Since the equivariant Euler class
of the normal bundle to a component of the fixed point set is never a zero-divisor,
we have $\eta\in\ker\phi_{\emptyset}$.
Then by Equation \ref{pols}, $$a+\eta\cdot (b_1+x)^{\sz}\in\mathcal{J}_S.$$
This completes the proof of Theorem \ref{eqcore}.
\end{eqcoreproof}

\vspace{-\baselineskip}
\begin{example}\label{projective}
For arbitrary $n$ and $\a$, suppose that $S$ is a maximal short subset.
Then Corollary \ref{ordcore} tells us that 
$H^*(U_S)\cong\Q[b_1]/\langle b_1^{n-2}\rangle$.
We conjecture that in this case we in fact have $U_S\cong \C P^{n-3}$.
\end{example}

\begin{example}\label{again}
Consider the core component pictured in Example \ref{ooh}.
By Theorem \ref{ordcore} and Remark \ref{integer}, 
$$H^*(U_S) \cong \Q[b_1,b_3,b_4,b_5]
\Bigg/
\left<\begin{array}{c}
b_3(b_1-b_3),\hs b_4(b_1-b_4),\hs b_5(b_1-b_5),\hs b_3b_4,\hs b_3b_5,\hs b_4b_5,\\
b_1(b_1-b_3-b_4),\hs  b_1(b_1-b_3-b_5),\hs b_1(b_1-b_4-b_5)
\end{array}\right>,$$
where $b_1$ is the fundamental class of $M_S$,
and $b_3$, $b_4$, and $b_5$ are the negatives of the fundamental classes of the curves
labeled $123$, $124$, and $125$, respectively.
Because the transverse intersection of two complex varieties is positive,
we know that $-b_1b_3[U_S] = 1$.
With respect to the basis $$\{b_1-b_3-b_4-b_5,\, b_3,\, b_4,\, b_5\},$$
the intersection form on $H^2(U_S)$ is represented by the matrix
$$\left(\begin{array}{cccc}
1&&&\\
&-1&&\\
&&-1&\\
&&&-1\end{array}\right).$$
Hence $U_S$ is homeomorphic to the blow-up of $\C P^2$ at three points.
\end{example}

\begin{example}
Using the same $\a=0\oplus (1,1,3,3,3)$,
consider the short subset $S=\{1,3\}$.
In this case, Theorem \ref{ordcore} tells us that
$$H^*(U_S)\cong\Q[b_1,b_2]/\left<b_1^2, b_2(b_1-b_2)\right>.$$
With respect to the basis $\{b_1-b_2,\, b_2\}$,
the intersection form on $H^2(U_S)$ is represented by the matrix
$\footnotesize{\Big(\begin{array}{cc}
-1&0\\
0& 1
\end{array}
\Big)}$,
hence $U_S$ is homeomorphic to the blow-up of $\C P^2$ at a single point.
\end{example}
\end{section}

\footnotesize{

}

\end{spacing}

\end{document}